\input amstex
\documentstyle{amsppt}
\magnification=\magstep1                        
\hsize6.5truein\vsize8.9truein                  
\NoRunningHeads
\loadeusm

\magnification=\magstep1                        
\hsize6.5truein\vsize8.9truein                  
\NoRunningHeads
\loadeusm

\document
\topmatter

\title
The sharp Remez-type inequality for even trigonometric polynomials on the period
\endtitle

\author Tam\'as Erd\'elyi
\endauthor

\address Department of Mathematics, Texas A\&M University,
College Station, Texas 77843, College Station, Texas 77843 \endaddress

\thanks {{\it 2010 Mathematics Subject Classifications.} 11C08, 41A17, 26C10, 30C15}
\endthanks

\keywords
trigonometric polynomials, Remez-type inequalities, geometry of polynomials
\endkeywords

\date September 7, 2018
\enddate

\email terdelyi\@math.tamu.edu
\endemail

\abstract
We prove that 
$$\max_{t \in [-\pi,\pi]}{|Q(t)|} \leq T_{2n}(\sec(s/4)) = 
\frac 12 ((\sec(s/4) + \tan(s/4))^{2n} + (\sec(s/4) - \tan(s/4))^{2n})$$
for every even trigonometric polynomial $Q$ of degree at most $n$ with 
complex coefficients satisfying
$$m(\{t \in [-\pi,\pi]: |Q(t)| \leq 1\}) \geq 2\pi-s\,, \qquad s \in (0,2\pi)\,,$$
where $m(A)$ denotes the Lebesgue measure of a measurable set $A \subset {\Bbb R}$ and 
$T_{2n}$ is the Chebysev polynomial of degree $2n$ on $[-1,1]$ defined by $T_{2n}(\cos t) = \cos(2nt)$
for $t \in {\Bbb R}$.
This inequality is sharp. We also prove that 
$$\max_{t \in [-\pi,\pi]}{|Q(t)|} \leq T_{2n}(\sec(s/2)) = 
\frac 12 ((\sec(s/2) + \tan(s/2))^{2n} + (\sec(s/2) - \tan(s/2))^{2n})$$
for every trigonometric polynomial $Q$ of degree at most $n$ with complex coefficients 
satisfying $$m(\{t \in [-\pi,\pi]: |Q(t)| \leq 1\}) \geq 2\pi-s\,, \qquad s \in (0,\pi)\,.$$
\endabstract

\endtopmatter

\head 0. Foreword \endhead

I started my Ph.D. program in June, 1987, at The Ohio State University, as a student of Paul Nevai. 
In August, 1987, Paul Nevai moved from Columbus, Ohio, to Columbia, South Carolina, to spend the 
school year 1987--88 at the University of South Carolina, and each of his students followed him. 
This is how I met Mr. (Yingkang) Hu, who was a student of Ron DeVore at that time. It was a 
somewhat short but very exciting time for me at the University of South Carolina. Not only had I 
passed each of my necessary exams in my Ph.D. program but I had the possibility to learn from and interact 
with some of the very best researchers in approximation theory. My office neighbors were Paul Nevai, 
Ron DeVore, and George Lorenz. George Lorentz happened to be Ron DeVore's office neighbor on one hand 
and my office neighbor on the other hand to finalize the book entitled ``Constructive Approximation" 
he was writing jointly with Ron DeVore. It was a privilege for me that the authors of this 
fundamental book in approximation theory asked me to read some sections of their book in preparation 
and took my remarks seriously. I was also lucky to have talented, motivated, and helpful 
class-mates such as Yingkang Hu, George Kyriazis, Jiaxiang (John) Zhang, and Shing-Whu Jha around. 
During the school year 1987--88 Yingkang was a good friend of mine. We spent much time together as 
class-mates as well as friends. I liked to eat Chinese food in the company of Mr. Hu, Jiaxiang Zhang, 
Shing-Whu Jha, and some of my other Chinese friends. We all took Paul Nevai's class on 
Orthogonal Polynomials. We traveled quite much together in this school year to attend various
conferences, but Yinkang, Jiaxiang, and me traveled together as tourists as well in the spring break.   
In August, 1988, Paul Nevai moved back to Columbus, and I followed him to write my dissertation at 
The Ohio State University. I defended my Ph.D. dissertation at the University of South Carolina in 
May, 1989. A central part of my dissertation was a Remez-type inequality for (generalized) algebraic 
and trigonometric polynomials and their applications in the proofs of various other inequalities for 
generalized polynomials. I was not able to prove the sharp Remez-type inequality for trigonometric 
polynomials, but I have managed to prove the sharp upper bound for $|Q(0)|$  for even trigonometric 
polynomials $Q$ satisfying 
$$m\{t \in (0,2\pi]:|Q(t)| \leq 1\}) \geq 2\pi - s\,, \qquad s \in (0,2\pi)\,,$$ 
where $m(A)$ denotes the Lebesgue measure of a set $A \subset {\Bbb R}$. After some 30 years 
I have revisited the topic of my dissertation. While the sharp Remez-type inequality for 
trigonometric polynomials remains open, I am glad to prove the sharp Remez-type inequality at least 
for even trigonometric polynomials in this volume dedicated to the memory of Dr. Hu. After 
defending my thesis I met Yingkang very rarely, but we were keeping in touch. In March, 2002, 
Yingkang visited us at Texas A\&M University, and gave a talk in our Center of Approximation Theory 
Seminar with the title ``Equivalence of Moduli of Smoothness", a topic on which he published two 
papers in Journal of Approximation Theory. During his visit we watched the Academy Award Ceremony 
in my apartment as he was very much interested in some of the latest movies. This may have been the 
last time I saw him and I was very sad to learn about his death in March, 2016.

\head 1. Introduction \endhead
Let ${\Cal T}_n$ denote the set of all real trigonometric polynomials of degree at most $n$. 
Let ${\Cal T}_n^c$ denote the set of all complex trigonometric polynomials of degree at most $n$.
Let $K := {\Bbb R} \enskip(\text {mod}\,\, 2\pi)$. Let $m(A)$ denote the Lebesgue measure of 
a measurable set $A \subset {\Bbb R}$. For $s \in (0,2\pi)$ let 
$${\Cal T}_n(s) := \{Q \in {\Cal T}_n: m(\{t \in K: |Q(t)| \leq 1\}) \geq 2\pi-s\}$$
and
$${\Cal T}_n^c(s) := \{Q \in {\Cal T}_n^c: m(\{t \in K: |Q(t)| \leq 1\}) \geq 2\pi-s\}\,.$$
Let ${\Cal P}_n$ denote the set of all algebraic polynomials of degree at most $n$ with real coefficients. 
Let ${\Cal P}_n^c$ denote the set of all algebraic polynomials of degree at most $n$ with 
complex coefficients.
Let $T_n$ be the Chebysev polynomial of degree $n$ on $[-1,1]$ defined by $T_n(\cos t) = \cos(nt)$ 
for $t \in K$. 
For real numbers $a < b$ and a Lebesgue measurable set $A \subset [a,b]$ let
$$\mu_{[a,b]}(A) := \int_A{\frac{(b-a)/2}{\sqrt{((b-a)/2)^2 - (x-(a+b)/2)^2}} \, dx}\,.$$

The classical Remez inequality states that if $P \in {\Cal P}_n$, $s \in (0,2)$, and
$$m \left( \left\{x \in [-1,1]: |P(x)| \leq 1 \right\} \right) \geq 2-s\,,$$
then
$$\max_{x \in [-1,1]}{|P(x)|} \leq T_n \left(\frac{2+s}{2-s} \right)\,,$$
where $T_n$ defined by
$$T_n(x) := \cos (n \arccos x)\,, \qquad x \in [-1,1]\,,$$
is the Chebyshev polynomial of degree $n$.
This inequality is sharp and
$$ T_n \left(\frac{2+s}{2-s} \right) \leq \exp(\min\{5ns^{1/2},2n^2s\})\,, \qquad s \in (0,1]\,.$$
Remez-type inequalities for various classes of functions have been studied by several authors, and 
they have turned out to be applicable and connected to various problems in approximation theory. 
See [1]--[34], for example.

In [16] we proved that 
$$|Q(0)| \leq T_{2n}(\sec(s/4)) \tag 1.1$$
for every even $Q \in {\Cal T}_n^c(s)$ and $s \in (0,2\pi)$. However, 
$0$ is a special point in the study of even trigonometric polynomials $Q \in {\Cal T}_n^c(s)$, so the 
question whether or not 
$$\max_{t \in K}{|Q(t)|} \leq T_{2n}(\sec(s/4)) \tag 1.2$$
holds at least for all even trigonometric polynomials $Q \in {\Cal T}_n^c(s)$ remained open, 
while it was speculated that (1.2) may hold for all $Q \in {\Cal T}_n^c(s)$.   
In this paper we show that (1.2) holds for all even $Q \in {\Cal T}_n^c(s)$, while it remains open whether 
or not (1.2) holds for all $Q \in {\Cal T}_n^c(s)$.

\head 2. New Results \endhead

Our first result is a sharp Remez-type inequality for even trigonometric polynomials  
with complex coefficients. 

\proclaim{Theorem 2.1} 
We have
$$\max_{t \in K}{|Q(t)|} \leq T_{2n}(\sec(s/4)) = 
\frac 12 ((\sec(s/4) + \tan(s/4))^{2n} + (\sec(s/4) - \tan(s/4))^{2n})$$
for every even $Q \in {\Cal T}_n^c(s)$ and $s \in (0,2\pi)$. Equality holds if and only if 
$$\{t \in K: |Q(t)| \leq 1\} = [s/2, 2\pi-s/2]$$ 
and $Q \in {\Cal T}_n(s)$ is of the form
$$Q(t) = T_{2n}\left(\frac{\cos(t/2)}{\cos(s/4)}\right)\,.$$ 
\endproclaim

Our next result is a Remez-type inequality for odd trigonometric polynomials with complex coefficients.

\proclaim{Theorem 2.2}
We have
$$\max_{t \in K}{|Q(t)|} \leq  
\frac 12 ((\sec(s/4) + \tan(s/4))^{2n} + (\sec(s/4) - \tan(s/4))^{2n}) + \frac{1}{\sqrt 2}$$
for every odd $Q \in {\Cal T}_n^c(s)$ and $s \in (0,2\pi)$. 
\endproclaim

Theorem 2.1 implies he following result for all trigonometric polynomials with complex coefficients.

\proclaim{Theorem 2.3}
We have
$$\max_{t \in K}{|R(t)|} \leq T_{2n}(\sec(s/2)) = 
\frac 12 ((\sec(s/2) + \tan(s/2))^{2n} + (\sec(s/2) - \tan(s/2))^{2n})$$
for every $R \in {\Cal T}_n^c(s)$ and $s \in (0,\pi)$. 
\endproclaim

Note that Theorem 2.2 was proved in [16] only for all $R \in {\Cal T}_n(s)$ 
(rather than $R \in {\Cal T}_n^c(s)$) and $s \in (0,\pi)$. It remains open whether or not 
Theorem 2.1 can be extended to all $Q \in {\Cal T}_n^c(s)$ and $s \in (0,2\pi)$.

\head 3. Lemmas \endhead 

Our first lemma is the conclusion (13) of Section 7 in [16]. This deep result plays a central 
role in the proof of Theorem 2.1.  

\proclaim{Lemma 3.1} 
We have
$$|U(0)| \leq T_{2n}(\sec(s/4))$$
for every even $U \in {\Cal T}_n(s)$ and $s \in (0,2\pi)$. Equality holds if and only if
$$\{t \in K: |U(t)| \leq 1\} = [s/2,2\pi-s/2]$$
and
$$U(t) = T_{2n}\left(\frac{\cos(t/2)}{\cos(s/4)}\right)\,.$$
\endproclaim

Our next lemma extends Lemma 3.1 to all even trigonometric polynomials with complex coefficients. 

\proclaim{Lemma 3.2}
We have
$$|Q(0)| \leq T_{2n}(\sec(s/4))$$
for every even $Q \in {\Cal T}_n^c(s)$ and $s \in (0,2\pi)$. Equality holds if and only if
$$\{t \in K: |Q(t)| \leq 1\} = [s/2,2\pi-s/2]$$
and
$$Q(t) = T_{2n}\left(\frac{\cos(t/2)}{\cos(s/4)}\right)\,.$$
\endproclaim

\demo{Proof}
Let $Q \in {\Cal T}_n^c(s)$. Choose $c \in {\Bbb C}$ with $|c|=1$ so that $cQ(0)$ is real.
Define $U \in {\Cal T}_n(s)$ by
$$U(t) := \text {\rm Re}(cQ(t))\,, \qquad t \in K\,.$$
Applying Lemma 3.1 to $U \in {\Cal T}_n(s)$ we obtain
$$|Q(0)| = |S(0)| \leq T_n(\sec(s/4))\,,$$
and equality holds if and only if
$$\{t \in K: |Q(t)| \leq 1\} = [-\pi+s/2,\pi-s/2]$$
and
$$Q(t) = T_{2n}\left(\frac{\cos(t/2)}{\cos(s/4)}\right)\,.$$
\qed \enddemo

Observe that every even $U \in {\Cal T}_n^c$ can be written as $U(t) = P(\cos t)$, where $P \in {\Cal P}_n^c$. 
In terms of algebraic polynomials Lemma 3.1 can be formulated as follows.

\proclaim{Lemma 3.3}
We have
$$|P(1)| \leq T_{2n}(\sec(s/4))$$
for every $P\in {\Cal P}_n^c$ such that 
$$\mu_{[-1,1]}(\{x \in [-1,1]: |P(x)| \leq 1\}) \geq \pi-s/2\,.$$
Equality holds if and only if
$$\{x \in [-1,1]: |P(x)| \leq 1\} = [-1,\cos(s/2)]$$
and
$$P(\cos t) = T_{2n}\left(\frac{\cos(t/2)}{\cos(s/4)}\right)\,.$$
\endproclaim

Associated with $P\in {\Cal P}_n^c$ let ${\widetilde P} \in {\Cal P}_n^c$ be defined by 
${\widetilde P}(x) = P(-x)$. Observing that 
$$\mu_{[-1,1]}(\{x \in [-1,1]: |P(x)| \leq 1\}) = \mu_{[-1,1]}(\{x \in [-1,1]: |{\widetilde P}(x)| \leq 1\})$$ 
we obtain the following.

\proclaim{Lemma 3.3*}
We have
$$|P(-1)| \leq T_{2n}(\sec(s/4))$$
for every $P\in {\Cal P}_n^c$ such that
$$\mu_{[-1,1]}(\{x \in [-1,1]: |P(x)| \leq 1\}) \geq \pi-s/2\,.$$
Equality holds if and only if
$$\{x \in [-1,1]: |P(x)| \leq 1\} = [-\cos(s/2),1]$$
and
$$P(-\cos t) = T_{2n}\left(\frac{\cos(t/2)}{\cos(s/4)}\right)\,.$$
\endproclaim

We need to transform the above two lemmas linearly from the interval $[-1,1]$ to the interval $[a,b]$. 

\proclaim{Lemma 3.4}
We have
$$|P(b)| \leq T_{2n}(\sec(s/4))$$
for every $P\in {\Cal P}_n^c$ such that
$$\mu_{[a,b]}(\{x \in [a,b]: |P(x)| \leq 1\}) \geq \frac{b-a}{2}\,(\pi-s/2)\,.$$
Equality holds if and only if
$$\{x \in [a,b]:|P(x)| \leq 1\} = [a,(a+b)/2+((b-a)/2)\cos(s/2)]$$
and
$$P(((b-a)/2)\cos t + (a+b)/2) = T_{2n}\left( \frac{\cos(t/2)}{\cos(s/4)}\right)\,.$$
\endproclaim

\proclaim{Lemma 3.4*}
We have
$$|P(a)| \leq T_{2n}(\sec(s/4))$$
for every $P\in {\Cal P}_n^c$ such that
$$\mu_{[a,b]}(\{x \in [a,b]: |P(x)| \leq 1\}) \geq \frac{b-a}{2}\,(\pi-s/2)\,.$$
Equality holds if and only if
$$\{x \in [a,b]: |P(x)| \leq 1\} = [(a+b)/2-((b-a)/2)\cos(s/2),b]$$
and
$$P((-(b-a)/2)\cos t + (a+b)/2) = T_{2n}\left( \frac{\cos(t/2)}{\cos(s/4)}\right)\,.$$
\endproclaim

Observing that $P \in {\Cal P}_n^c$ implies ${\widetilde P} \in {\Cal P}_n^c$, where 
${\widetilde P}(x) = P(-x)$ we obtain the following. 

\proclaim{Lemma 3.5}
We have
$$\mu_{[1-2r,1]}(A) > r^{1/2} \mu_{[-1,1]}(A)$$
for every $r \in (0,1)$ and for every Lebesgue measurable set $A \subset [1-2r,1]$ with $m(A) > 0$.
\endproclaim

\demo{Proof} 
Let $r \in (0,1)$ and let $A \subset [1-2r,1]$ be a Lebesgue measurable set with $m(A) > 0$. 
We have 
$$\split \frac{\mu_{[1-2r,1]}(A)}{\mu_{[-1,1]}(A)} & = 
\frac{\displaystyle{\int_A{\frac{r}{\sqrt{r^2 - (x-(1-r))^2}} \, dx}}}{\displaystyle{\int_A{\frac{1}{\sqrt{(1-x^2}} \, dx}}} 
= \frac{ \displaystyle{\int_A{\frac{r}{\sqrt{(1-x^2)}} \, \frac{\sqrt{(1-x^2}}{\sqrt{r^2 - (x-(1-r))^2}}} \, dx}}
{\displaystyle{\int_A{\frac{1}{\sqrt{1-x^2}}}\,dx}} \cr 
& > \frac{\displaystyle{\int_A{\frac{r}{\sqrt{1-x^2}}\,dx}} \, \min_{y \in [1-2r,1]}{\displaystyle{\frac{\sqrt{1-y^2}}{\sqrt{r^2 - (y-(1-r))^2}}}}}
{\displaystyle{\int_A{\frac{1}{\sqrt{1-x^2}} \, dx}}} \cr 
& \geq rr^{-1/2} = r^{1/2} \cr \endsplit$$
as
$$\min_{y \in [1-2r,1]}{\frac{\sqrt{1-y^2}}{\sqrt{r^2 - (y-(1-r))^2}}}
= \min_{y \in [1-2r,1]}{\frac{\sqrt{1+y}}{\sqrt{y - (1-2r)}}} = \frac{\sqrt{2}}{\sqrt{2r}} = r^{-1/2}\,.$$
\qed \enddemo

\proclaim{Lemma 3.6}
We have
$$\mu_{[-1,1-2r]}(A) > (1-r)^{1/2} \mu_{[-1,1]}(A)$$
for every $r \in (0,1)$ and for every Lebesgue measurable set $A \subset [-1,-1+2r]$ with $m(A) > 0$. 
\endproclaim

\demo{Proof}
Let $r \in (0,1)$ and let $A \subset [-1,1-2r]$ be a Lebesgue measurable set with $m(A) > 0$.. 
Applying Lemma 3.5 with $A$ replaced by   
$-A := \{-x:x \in A\} \subset [1-2(1-r),1]$ and $r \in (0,1)$ replaced by $1-r \in (0,1)$, we get   
$$\mu_{[-1,-1+2r]}(A) = \mu_{[1-2(1-r),1]}(-A) > (1-r)^{1/2} \mu_{[-1,1]}(-A) 
= (1-r)^{1/2} \mu_{[-1,1]}(A)\,.$$
\qed \enddemo

\proclaim{Lemma 3.7}
We have
$$\arccos(1-2r) > \pi r^{1/2} \qquad and \qquad \pi - \arccos(1-2r) > \pi (1-r)^{1/2}$$
for all $r \in (0,1)$. 
\endproclaim

\demo{Proof}
The first inequality is equivalent to $\sin t > (2/\pi)t$ for all $t = (\pi/2)r^{1/2} \in (0,\pi/2)$.
The second inequality follows from the first one as
$$\pi - \arccos(1-2r) = \arccos(-1+2r) = \arccos(1-2(1-r)) > \pi (1-r)^{1/2}$$
for all $r \in (0,1)$.
\qed \enddemo  

\head 4. Proof of Theorem 2.1 \endhead

\demo{Proof of Theorem 2.1}
Let $Q \in {\Cal T}_n^c(s)$ be even with $s \in (0,2\pi)$. Then $Q$ is of the form 
$Q(t) = P(\cos t)$, where $P \in {\Cal P}_n^c$ satisfies
$$\mu_{[-1,1]}(\{x \in [-1,1]: |P(x)| \leq 1\}) \geq \pi-s/2\,. \tag 4.1$$
We want to prove that
$$|P(\alpha)| \leq T_{2n}(\sec(s/4)) \tag 4.2$$
for all $\alpha \in [-1,1]$. 

If $\alpha=1$ then (4.2) holds by Lemma 3.3. If $\alpha=-1$ then (4.2) holds by Lemma 3.3*. 
Assume now that $\alpha = 1-2r \in (-1,1)$, that is, $\alpha = 1-2r$ with $r \in (0,1)$. It follows from (4.1)  
and the definition of the measure $\mu_{[-1,1]}$ that
$$m(\{t \in [0,\pi]:|P(\cos t)| \leq 1)\}) \geq \pi-s/2\,,$$
and hence 
we have either 
$$m(\{t \in [0,\arccos(1-2r)]:|P(\cos t)| \leq 1)\}) \geq \frac{\pi-\arccos(1-2r)}{\pi}(\pi-s/2) \tag 4.3$$
or
$$m(\{t \in [\arccos(1-2r),\pi]:|P(\cos t)| \leq 1)\}) \geq \frac{\arccos(1-2r)}{\pi}(\pi-s/2)\,. \tag 4.4$$  
For the sake of brevity let 
$$A_1 := \{x \in [-1,1-2r]: |P(x)| \leq 1\} \quad \text {\rm and} \quad 
A_2 := \{x \in [1-2r,1]:|P(x)| \leq 1\}\,. \tag 4.5$$
By (4.3), (4.4), and (4.5) we have either
$$\mu_{[-1,1]}(A_1) \geq \frac{\pi-\arccos(1-2r)}{\pi}(\pi-s/2) \tag 4.6$$
or 
$$\mu_{[-1,1]}(A_2) \geq \frac{\arccos(1-2r)}{\pi}(\pi-s/2)\tag 4.7$$

Suppose that (4.6) holds. Then, combining Lemmas 3.6 and 3.7 we obtain  
$$\split \mu_{[-1,1-2r]}(A_1) > & (1-r)^{1/2}\mu_{[-1,1]}(A_1) \geq (1-r)^{1/2}\frac{\pi-\arccos(1-2r)}{\pi}(\pi-s/2) \cr 
> & (1-r)^{1/2}(1-r)^{1/2}(\pi-s/2) \geq (1-r)(\pi-s/2)\,, \cr \endsplit$$
and 
$$|P(\alpha)| = |P(1-2r)| < T_{2n}(\sec(s/4))$$
follows from Lemma 3.4.

Suppose that (4.7) holds. Then, combining Lemmas 3.5 and 3.7 we obtain
$$\split \mu_{[1-2r,1]}(A_2) > & \, r^{1/2}\mu_{[-1,1]}(A_2) \geq r^{1/2}\frac{\arccos(1-2r)}{\pi}(\pi-s/2) \cr 
> \, & r^{1/2}r^{1/2}(\pi-s/2) \geq r(\pi-s/2)\,, \cr \endsplit$$
and
$$|P(\alpha)| = |P(1-2r)| < T_{2n}(\sec(s/4))$$
follows from Lemma 3.4*.
\qed \enddemo

\demo{Proof of Theorem 2.2}
Let $Q \in {\Cal T}_n^c(s)$ be odd. Then $R$ defined by $R(t) = 2|Q(t)|^2-1$ is even and 
$R \in {\Cal T}_{2n}^c(s)$. Applying Theorem 2.1 to $R$ we obtain
$$\max_{t \in K}{|2|Q(t)|^2-1|} = \max_{t \in K}{|R(t)|} \leq
\frac 12 ((\sec(s/4) + \tan(s/4))^{4n} + (\sec(s/4) - \tan(s/4))^{4n})\,,$$
and hence
$$\split \max_{t \in K}{|Q(t)|} \leq &  
\left( \frac 14 ((\sec(s/4) + \tan(s/4))^{4n} + (\sec(s/4) - \tan(s/4))^{4n}) + \frac 12 \right)^{1/2} \cr
\leq & \frac 12 ((\sec(s/4) + \tan(s/4))^{2n} + (\sec(s/4) - \tan(s/4))^{2n}) + \frac{1}{\sqrt 2} \cr \endsplit$$
\qed \enddemo

\demo{Proof of Theorem 2.3}
Observe that if $R \in {\Cal T}_n^c(s)$, $s \in (0,\pi)$, and $Q$ is defined by $Q(t) = \frac 12(R(t) + R(-t))$,
then $R \in {\Cal T}_n^c(2s)$. Hence Theorem 2.1 implies that
$$|R(0)| = |Q(0)| \leq T_{2n}(\sec(s/2)) \tag 4.8$$
holds for every even $R \in {\Cal T}_n^c(s)$ and $s \in (0,\pi)$. The theorem now follows by a simple shift. 
Namely, observe that $R \in {\Cal T}_n^c(s)$ implies that $R_a \in {\Cal T}_n^c(s)$ where $R_a$ is defined by 
$R_a(t) := R(t+a)$ for $t \in K$ and $a \in K$. Hence (4.8) implies that
$$|R(a)| = |R_a(0)| \leq T_{2n}(\sec(s/2))$$  
holds for every even $R \in {\Cal T}_n^c(s)$, $s \in (0,\pi)$, and $a \in K$.
\qed \enddemo

\enddocument